\documentstyle{amsppt}
\headline={\hfil}
\footline={\ifnum\pageno<2 \hfil \else \hss\folio\hss \fi}
\output={\plainoutput}
%\printoptions
\magnification=\magstep1

\hsize=6.5truein
\vsize=9truein
\abovedisplayskip=15pt plus 3pt minus 3pt
\abovedisplayshortskip=10pt plus 3pt
\belowdisplayskip=15pt plus 3pt minus 3pt
\belowdisplayshortskip=10pt plus 3pt minus 2pt
\NoBlackBoxes
\def\qed{\hfill $\square$}

\parindent = 20pt      

\textfont\ttfam=\tentt %

 % \q abreviation for \quad

\def\noi{\noindent}

\def\g{\frak g}

\def\fX{ X}
\def\X{X}

\def\Y{ Y}

\def\C{L}
\def\fg{\frak g}

\def\A{\Bbb A^1}

\def\fl{\frak l}

\def\Ralg{R \text{-alg}}
\def\s{\sigma}

\def\AdG{\text{\bf Ad}(G)}
\def\AdB{\text{\bf Ad}(B)}
\def\AdT{\text{\bf Ad}(T)}
\def\Aut{\text{\bf Aut}}
\def\Out{\text{\bf Out}}
\def\AutA{\text{\bf Aut}(A)}
\def\AutG{\text{\bf Aut}(G)}
\def\OutG{\text{\bf Out}(G)}

\def\Spec{\text{\rm Spec}\,}
\def\Gal{\text{\rm Gal}}

\def\Pic{\text{\rm Pic}\,}
\def\Aut{\text{\rm \bf Aut}}
\def\Out{\text{\rm \bf Out}}
\def\rmAut{\text{\rm  Aut}}
\def\rmOut{\text{\rm  Out}}
\def\0{\{0\}}
\def\et{\text{\'et}}

\topmatter  
\title 
\centerline{Affine Kac-Moody Lie algebras as torsors }
\centerline{over the punctured line}
\endtitle

\author
A. Pianzola\footnotemark"$^1$" 
 \endauthor

\abstract\nofrills{Abstract.}
We interpret and develop a theory of loop algebras  as torsors (principal homogeneous 
spaces)  over $\Spec(k[t ,t^{-1}])$ . As an application, we recover 
Kac's realization of affine Kac-Moody Lie algebras.
\endabstract

$\footnotetext"$^1$"{Supported by NSERC operating grant A9343.}$

\affil{\centerline{Department of  Mathematical Sciences}
\centerline{ University of Alberta}
\centerline  {Edmonton, Alberta,  T6G 2G1}
\centerline {Canada}}
\endaffil


\endtopmatter

\noi {\bf Introduction.} There is a beautiful construction of Victor Kac's, realizing affine 
Kac-Moody Lie algebras over the complex numbers as (twisted) loop algebras. The construction 
gives explicit 
 generators for the algebras,
which are then shown to satisfy the relations corresponding to the affine Cartan 
matrix at hand.

In this short note, we propose to look at loop algebras in a 
completely different way. The basic idea is to view loop algebras  as 
algebras over a ring of Laurent polynomials $R$, all of which 
become isomorphic after a flat covering  $R \rightarrow S.$ Thus loop 
algebras become torsors over $\Spec (R)$ under the group of 
automorphisms of the algebra at hand. Since this 
point of view applies to
arbitrary algebras and base fields, we are able to obtain some rather general new 
results about loop algebras in {\bf 
8} and {\bf 10}. As an application, we show how to recover from this  Kac's 
original result. This is done in {\bf 11.}

We begin with a review of loop algebras in {\bf 1}, and then recall 
Kac's construction in {\bf 2}. This is followed by some abstract results on algebraic 
groups ({\bf 3} through {\bf 5}) which are later needed. The 
description of loop algebras as torsors is given in {\bf 7}. 
\bigskip
\noi {\bf 0. Conventions and notation.} Throughout  this note $k$ will denote a field. If $G$ is a $k$-group and 
 $\X = \Spec(R)$   
 an affine $k$-scheme, the 
 $\X$-group $G \times_{\Spec(k)} \X$ will be denoted by $G_{\X}$ or 
 $G_{R}$ .  For an $X$-group $F$ the C\v ech cohomology $\check 
 H^{1}(X_{\et}, F)$ on the \'etale site of $X$ wil be denoted simply by $
 H^{1}(X, F)$. For terminology 
and results about schemes and principal homogeneous spaces (torsors), the reader is referred to 
[DG] Ch. 3.4, [Mln] Ch. 3, and [SGA1].
\medskip
We will throughout work with two copies of $k$-algebras of Laurent polynomials that for 
convenience we denote by 
$$R = k[t , t^{-1}] \,\,\text{and}\,\, S = k[z , z^{-1}].$$

For each  positive 
integer $m$ we view the ring $S$ as an $R$-algebra via $t \mapsto z^{m}$. This 
algebra will be denoted by $S_{m}$. As usual in what follows $\, \bar{} 
: \Bbb Z \rightarrow \Bbb Z/m\Bbb Z$ will denote the canonical map.

\medskip
\noi {\bf 1. Loop algebras.} 
Let $A$ be an algebra over $k$ of a certain ``type'' (eg. 
associative, Lie, Jordan, etc) which we assumed comes equiped with a 
$\Bbb Z/m\Bbb Z$-grading $\Sigma.$ Thus 

$$A = \underset{0 \leq i < m}\to\oplus A_{\bar \imath }$$  
where the  $A_{\bar \imath }$ are subspaces of $A$ such that $ A_{\bar \imath_{1} } A_{\bar \imath_{2} } \subset A_{\bar 
\imath_{1} + \bar \imath_{2}}.$ We then define 
the {\it loop 
algebra} of $A$ with respect to the given grading $\Sigma$ by 

$$\C(A,\Sigma) := \underset{i \in \Bbb Z}\to\oplus A_{\bar \imath } \otimes 
z^i \subset A \otimes_{k}€ k[z , z^{-1}] := A \otimes_{k}S.$$

\noi Note  that loop algebras 
are in a natural way $k-$algebras of the same type as $A$ (eg.  
loop algebras of a Lie algebra  are  Lie algebras etc.), and are 
viewed as such in what follows.  If the grading is trivial, namely if $A_{\bar 0} = A$, then 
$\C(A,\Sigma) \simeq  A \otimes_{k}R$. Loop algebras isomorphic as 
$k$-algebras to $A 
\otimes_{k}R$ are said 
to be {\it trivial.}
\medskip
Suppose now that $k$ is algebraically closed of characteristic $0$, 
and fix  a compatible set $\pmb{\xi} = (\xi_{m})$  of  primitive  
roots of unity in $k$ (thus if $n = cm$ then ${\xi_{n}}^{c} = \xi_{m}$). Let $\s$ be an 
automorphism of $A$ of period $m$ (we  use periods to be able to compare  loop algebras coming from 
automorphisms of different order). 
Then  $A$ 
decomposes as the direct sum of  eigenspaces 
$A_{\overline \imath}$
where  $\s$ acts on $A_{\bar \imath}$ as scalar multiplication by 
${\xi_{m}}^i.$ We thus obtain a $\Bbb Z/m\Bbb Z$-grading $\Sigma$ of $A$ as 
above. Conversely, any such grading $\Sigma$ comes from a period $m$ 
automorphism $\sigma$ of $A.$ The  resulting loop algebra (which up to 
$R$-algebra  isomorphism is independent 
of the choice of period $m$ of $\s$)  will be 
denoted by $\C(A,\s, \pmb{\xi} ),$ or simply by $\C(\s)$ if $A$ and 
$\pmb{\xi}$ are fixed.   Note that $\C(\text{id}_{A})$ is trivial.
\medskip
Here is the most remarkable application of loop algebras. 
\proclaim{Theorem 2}    Let $\fg$ be a simple finite dimensional Lie 
algebra over $\Bbb C$, and let $^{-}$ 
denote the canonical map from  $\rmAut(\fg)$ onto the group 
$\rmOut(\fg)$ of outer automorphisms of $\g.$ Then.

\item{\rm (i)} $\C(\s) \simeq \C(\overline \s)$ as complex Lie 
algebras for all $\s \in 
\rmAut(\fg).$
\item{\rm (ii)} Let $\hat{\fl}$ be an affine Kac-Moody Lie algebra, and let $\fl$ be its derived
algebra modulo its centre. There exists  $\g$  and   $\pi 
\in \rmOut(\fg)$ as above  
 such that
$\fl \simeq \C(\pi).$
\item{\rm (iii)} In \text{ \rm (ii)} above, $\g$ is unique up to 
isomorphism, and $\pi$ unique up to conjugacy in  $\rmOut(\fg)$.

\endproclaim
\demo{Proof} Parts (i) and (ii) are due to Kac (see [Kac] Theorem 
8.3]). For (iii) one needs the conjugacy theorem of Peterson and Kac 
[PK]
\qed
\enddemo

Loop algebras are thus concrete realizations of the affine 
Kac-Moody Lie algebras.  
We propose now to give
new insight into this theorem, as well as obtain new results about  
loop algebras in general, by interpreting such algebras as 
 torsors over the punctured line. To this end, we begin with some results on the cohomology of
  algebraic groups.
\medskip 
  The following result shows that the classical vanishing of $H^{1}$ 
  theorems of Steinberg and of Borel and Springer, hold for certain 
  semisimple group schemes over Dedekind domains.
  
\proclaim{Proposition 3} Let $D$ be a Dedekind domain and  $K$  
its field of quotients. Set $\X = 
\Spec(D).$ Let $G$ be a quasisplit semisimple connected $\X$-group, $B$  a Borel 
subgroup of $G$, and $T$ a maximal torus of $B$.

\item{\rm (i)}  The natural map 
$H^{1}(\X, \, T) \rightarrow H^{1}(\X, \, G)$ induces a surjection 
$H^{1}(\X, \, T) \rightarrow \, \text{\rm Ker}\,(H^{1}(\X, \, G) \rightarrow 
H^{1}(K, \, G_{K}€)).$ 
\noi  In particular,  this kernel is trivial whenever $H^{1}(\X, \, T) $ is 
trivial.

\item{\rm (ii)} Assume $G$ is semisimple of either adjoint or simply connected 
type. If all connected \'etale coverings of $\X$ have trivial Picard 
group, then $H^{1}(\X, \, T) $ is trivial.

\item{\rm (iii)} Assume $G$ and $\X$ are as in (ii). If $K$ is of cohomological 
dimension 1, then $H^{1}(\X, \, G)$ is trivial.

\item{\rm (iv)} Assume that $\X$ and $K$ are as in (iii). If $G$ is 
split
semisimple then the canonical map 
\newline $H^{1}(\X, \, \Aut(G)) \rightarrow H^{1}(\X, \,\Out(G))$ is bijective.

\endproclaim
\demo{Proof} (i) This is the extent of Satz 3.2 in 
[Hrd]. Here is  another 
proof based on an idea (used in [CTO] in 
the case of a base field) more in tune with the spirit of this note. A torsor $\Y$ on the kernel in question 
is rationally trivial, i.e. 
it admits a section over a non-empty Zariski open of $\X.$ Consider the 
exact sequence. 
$$ 1 \rightarrow B \rightarrow G \rightarrow G/B \rightarrow 1$$
as well as the contracted product $\Y \times^{G} G/B.$ It is clear 
that the structure morphism $\kappa:\Y \times^{G} G/B \rightarrow \X$ admits 
a rational section. Since $\kappa$ is proper and $\X$ is one 
dimensional and regular, this section  extends to all of $\X$ ([EGA] 
Ch. II Cor. 7.3.6). By 
[DG] Ch.III 4.4.6 the $G$-torsor $\Y$ comes from $B,$ thence from 
$T$ given that
$H^{1}(\X, \, \text{rad}^{\text{\rm u}}(B)) = \0$ ([SGA3] XXII Cor. 
5.9.7). Note that this proof also works in the reductive case. That 
the elements of $H^{1}(X, \, T)$ are rationally trivial is prooved in 
section 1.4 of [Hrd]. 

(ii)  By [SGA3] XXIV 3.13-3.15 and [Hrd] 1.4
there exists
a finite family of connected \'etale coverings  $\X_i/\X$  of $\X$   such that
$$
 T = \prod_i R_{\X_i/\X} ({\Bbb G}_m),
$$
where the $ R_{\X_i/\X}$ are Weil restrictions. By Shapiro's lemma  one has
$$
 H^1(\X,\, T) = \prod_i  H^1(\X_{i}, {\Bbb G}_m)=
\prod_i \Pic(\X_i)=\0.
$$

(iii) By classical results  of Steinberg and of Borel and Springer  $H^{1}(K, \, G_{K}€)$ is 
trivial. 
(See [BS] 8.2, [Stb], and [JPS1] 
Ch. 3.1 and 3.2 ). Now (iii) follows from (i) and (ii)

(iv) Recall from [SGA3] XXV the existence of a split exact sequence of 
$\X$-groups
$$ 
1 \rightarrow \AdG  \rightarrow \Aut(G) \rightarrow
\Out(G) \rightarrow 1
$$
The group $\Out(G)$ is a finite  constant  
group, and 
admits a section $s : \Out(G) \rightarrow \Aut(G)$ whose
image consists of  elements of $\AutG$ that stabilize both $B$ and $T$ 
(We henceforth identify $\Out(G)$ with a subgroup of $\Aut (G)$ by 
means of this section). Passing to cohomology 
yields
$$H^{1}(\X, \, \AdG) \rightarrow 
H^{1}(\X, \, \Aut(G)) \rightarrow H^{1}(\X, \, \Out(G)).   $$
 Let $Y$ be an $X$-torsor under $\OutG$, and consider the groups 
 $_{Y}\AdG$ (resp.$ _{Y}\AdB$, $_{Y}\AdT$) obtained from $\AdG$ (resp. 
 $\AdB$, 
 $\AdT$) by twisting by $Y.$ Being a form of $\AdG$, the group 
 $_{Y}\AdG$ is  semisimple and of adjoint type as well. It is also 
 quasisplit by means of ($ _{Y}\AdB$, $_{Y}\AdT)$. By (iii) then,  
 $H^{1}(\X,  _{Y}\AdG)$ is trivial. From this it follows that  the map in question is 
injective. The surjectivity is clear because the existence of the 
section $s.$ \qed   

\medskip
\noi {\bf Note}. A priori $_{Y}\AdG$ is a sheaf of groups on $\X$. That it is an affine and 
smooth scheme over $X$ follows from descent. That its geometric fibers 
are reductive and connected follows from the analogous properties for 
$\AdG$. Thus $_{Y}\AdG$ {\it is} a reductive group in the sense of [SGA3].
Along similar lines $_{Y}\AdB$ is a Borel subgroup $\ldots$ 
\enddemo
\comment

To use [DG] and [Mln] one should work on the fppf site but since 
$_{Y}\AdG$ is smooth we can revert back to $\et$.  Of course it comes 
down to the existence of the quotient in Sch/X of $Y \times\AdG$ by 
the diagonal action of $\OutG$. By descent, this seems to be ok as long as $\OutG$ 
is affine (which it is, again by descent since $\OutG$ becomes 
constant after an \'etale covering). If in this generality all is correct, I 
would like to keep things as they are. 
\endcomment

\medskip
The following result shows that the  assumption on \'etale coverings 
made in part (ii) of Proposition 3, holds in a crucial case.

\proclaim{Proposition 4} {\rm (P. Gille)} \, Assume $k$ is 
 of characteristic $0$. Every  connected finite \'etale  covering of 
 $\Spec(k[t ,t^{-1}])$ has trivial Picard group. 
\endproclaim

\demo{ Proof}
Let $\Y \rightarrow \X := \Spec(k[t ,t^{-1}]€)$ be one such covering. Fix an element
$x \in X(k).$ From [SGA1]  Exp. IX Th. 6.3.1 together 
with a Theorem of Grauert-Remmert ([SGA4]  Exp. XI Th. 4.3), as well 
as  from 
[SGA1] Exp. XIII Cor. 2.12,  it follows  that the  fundamental
 group  $\Pi_{1}€(\X,x)$  is  the  semidirect product
 $(\text{inv lim}\,  \mu_n(k_s)) \rtimes \Gal(k_s/k)$,
where the $\mu_n$ come from the Kummer coverings.  
 There thus exists a positive integer $m$, a finite  Galois field extension $L/k$
containing a primitive $m$-th root of unity, and a subgroup $\Gamma$ of 
$\mu_m(L) \rtimes \Gal(L/k)$, 
such that $\Y = \Y_0 / \Gamma$ where 
$\Y_0$ is the $k$--variety defined by the morphism 
 $\Y_0= \X_{L} \rightarrow \Spec(L) \rightarrow
 \Spec(k)$. As the morphism $\Y_0 \rightarrow \Y$ is a Galois 
 covering, the beginning of
the Hochschild--Serre spectral sequence
$E_2^{p,q}=H^p(\Gamma, H^q({\Y_0}, {\Bbb G}_m)) \Longrightarrow
H^{p+q}(\Y, {\Bbb G}_m)$ yields an exact sequence
$$
0 \rightarrow H^1(\Gamma, H^0(\Y_0, {\Bbb G}_m)) \rightarrow  H^1(\Y,{\Bbb G}_m)
\rightarrow   H^1({\Y_0}, {\Bbb G}_m)^\Gamma.
$$

Since $H^1({\Y_0},{\Bbb G}_m)= \Pic(\Y_{0})= \0$,  we get an isomorphism
$H^1(\Gamma, H^0(\Y_0, {\Bbb G}_m)) \simeq \Pic(\Y)$.
One has an exact sequence of $\Gamma$--modules (and of
 $\mu_m(L) \rtimes \Gal(L/k)$--modules) 
$$
0 \rightarrow L^\times \rightarrow  H^0(\Y_0, {\Bbb G}_m) \rightarrow 
{\Bbb Z} \rightarrow 0,
$$
where ${\Bbb Z}$ has trivial $\Gamma$--action.
By Hilbert's theorem 90  one has $H^1(\Gamma, L^\times)= \0,$ and as 
$\Gamma$ is finite, one also has   $H^1(\Gamma, {\Bbb Z})= \0$. Thus  $H^1(\Gamma,
  H^0(\Y_0, {\Bbb G}_m))= \0$  and therefore $\Pic(\Y)= \0$ as desired.
\qed
\enddemo
Part (i) of the next result is an easy, but useful generalization of a result 
of [CTO].

\proclaim{Proposition 5} Let $k$ be an infinite perfect field and let $G$ be a 
smooth  connected linear algebraic $k$-group. Let $\X$ be a nonempty open subscheme 
of $\Spec(k[t]) = \A.$ Then.

\item{\text{\rm (i)}} The canonical map $H^{1}(\X, \, G_{\X}€) \rightarrow 
H^{1}(k(t), \, G_{k(t)})$ has trivial kernel. 
\item{\text{\rm (ii)}} 
$H^{1}(\X, G_{\X})$ is trivial in  the following three cases.

\item\item{\text{\rm (a)}} 
If $H^{1}(k(t), G_{k(t)})$ is trivial.

\item\item{\text{\rm (b)}} If $k$ is algebraically closed and of 
characteristic 0.

\item\item{\text{\rm (c)}} If $k$ is algebraically closed  and $G$ is reductive.

\endproclaim
\demo{Proof} Let us begin by showing that  every Zariski $G \text{-torsor}$ over $\X$ 
(i.e. a torsor that can be trivialized by a Zariski covering 
of $\X$) extends 
to  $\A$ (this much holds for any group scheme $G$). Indeed. The
underlying space  of $\X$ is obtained by removing a finite set $F$ of points 
from the affine line. Let now $\Y$ be a Zariski $G \text{-torsor}$ 
over $\X.$  Let $U$ be the intersection of a finite number of
non empty  open subschemes of $\X$ that cover  $\X,$ and  over which $\Y$ is 
trivial. Finally, let $Z$ be the trivial  $G \text{-torsor}$ over the open 
subscheme of  $\A$ corresponding to $U \cup F$. Then by gluing  
$\Y$  with $Z$ along $U$ 
we obtain a torsor over $\A$ as desired.

 (i) By reasoning as in  Th\'eor\`eme 2.1 of [CTO]  we may assume that 
$G$ is reductive and connected. A Theorem  of Nisnevich [Nsn] then shows that the kernel of the map in question  
is comprised precisely of Zariski $G$-torsors  over $\fX$. But we have seen that any such torsor extends to the full 
affine line. Now (i) follows from  [CTO] Corollaire 2.3 which asserts 
that (i) does hold for $\A$.

(ii) By the theorems of   Steinberg and  Borel-Springer 
mentioned above, $H^{1}(k(t), G_{k(t)})$ is trivial
 under the assumptions of 
either (b) 
or (c). Thus (b) and (c) reduce 
to (a), and this last holds by (i). 
  \qed

\medskip
\noi {\bf Remark 6.} Assume that $k$ contains a 
primitive $m$-th root of unity. Then the 
finite covering 
$\Spec(S_{m}¥) \rightarrow \Spec(R)$ is \'etale and in fact Galois, with 
Galois group $\Gamma \simeq \Bbb Z/m \Bbb Z$ (see {\bf 7} infra). Assume now that we are 
under one of the cases of Proposition 5 (ii). Then the usual non-abelian 
cohomology  $H^{1}(\Gamma, \, G(S_{m}¥))$ vanishes. If $k = \Bbb C$ and  $G$ 
is semisimple, we recover the observation made in [Kac] \S 8.9. 

  \enddemo

\medskip

\noi {\bf 7. Loop algebras as torsors.} We return now to our general setup of {\bf 
1} and consider an arbitrary $k$-algebra $A$, together with
$\Bbb Z/m\Bbb Z$-grading $\Sigma.$   For convenience we henceforth  denote
$\Spec(R) := \Spec(k[t,t^{-1}])$  by $\X$ and $\Spec(S) := 
\Spec(k[z,z^{-1}])$  by $\Y$

The loop algebra $ \C(\Sigma) = \C(A, \Sigma)$ is naturally an 
$R\text{-algebra}$, and it 
is not hard to show  that $\C(\Sigma) \otimes _{R} S_{m} \simeq 
A \otimes_{k}  S $ as $S$-algebras (see [ABP]). In other words, 
$\C(\Sigma)$ is an $S_{m}/R$-form 
of $A \otimes_{k}R.$ Since $S_{m}/R$ is  faithfully flat and finitely presented 
(fppf), $\C(\Sigma)$ is an $\Aut(A_{\X})$-torsor over $\X$. In this 
context, $\Aut(A_{X})$ stands for the sheaf of groups on the flat site 
of $\X$ that attaches 
to $\X'/\X$, the  group of $\Cal O(\X')$-algebra automorphism 
of $A \otimes_{k} \Cal O(X')$ of the same type as $A$.

The isomorphism class of the $R-$algebra $\C(\Sigma)$ is thus
an element of $H^{1}(\X_{\text{\rm fppf}}, \Aut(A_{\X}))$ that we denote by 
$\C^{1}€(\Sigma).$ To say that $\C^{1}€(\Sigma_{1}) = \C^{1}€(\Sigma_{2})$ 
is to say that 
$\C(\Sigma_{1})$ and $\C(\Sigma_{2})$ are
isomorphic as {\it algebras over} R. It is clear then that $\C(\Sigma_{1})$ and $\C(\Sigma_{2})$ are
a fortiori isomorphic as  algebras over $k.$ The converse does not 
hold in general, but as we shall see in {\bf 10} and {\bf 11} , it does hold in some 
very interesting cases.

Next we assume that the base field $k$ is algebraically closed and of 
characterist zero, and go on to describe explicitly how to construct 
$\C^{1}(\sigma)$ from $\s.$ 
The 
finite covering 
$Y \rightarrow X$ is \'etale and in fact Galois. Its 
Galois group $\Gamma$ can be identifyed with $\Bbb Z/m \Bbb Z$ where $1 + m\Bbb Z$ acts 
on $S$ via $t \rightarrow \xi t.$ This leads to a natural action 
of $\Gamma$ on the abstract group of 
automorphisms $\text{\rm Aut}A \otimes_{k} S := \Aut (A_{X})(Y)$ of the 
$S$-algebra $A \otimes_{k} S.$ The map $u: n + m\Bbb Z \mapsto 
\s^{-n}$ is clearly a $1$-cocycle in $Z^{1}(\Gamma,\Aut (A_{X})(Y)) = 
Z^{1}((Y/X)_{\text{\rm fppf}},\Aut (A_{X})),$ and a straightforward computation shows that the 
form corresponding to $u$ is precisely $\C(\sigma).$ To obtain then 
$\C^{1}(\sigma)$ we simply follow $u$ along the canonical maps 
$Z^{1}((Y/X)_{\text{\rm fppf}},\Aut (A_{X})) \rightarrow H^{1}((Y/X)_{\text{\rm fppf}},\Aut (A_{X})) \subset 
H^{1}(X_{\text{\rm fppf}},\Aut (A_{X}))$
\medskip
Kac's theorem states that for complex simple Lie algebras, the 
isomorphism class of a covering algebra $L(\s)$ depends only on the outer part 
of $\s.$ In addition, if $\s$ is inner then $L(\s)$ is trivial. The 
following two Propositions show this to be a particular instance of 
a rather general picture.
\medskip
\proclaim{Proposition 8 } Let $k$ be an infinite perfect field. Assume 
that the $k$-algebra $A$ is finite 
dimensional,    and that the 
linear algebraic $k$-group $\Aut(A)$ is smooth and connected. If 
either of the conditions  of Proposition 5 (ii) hold, 
then all loop algebras of $A$ are  trivial.

\endproclaim
\demo{Proof} By  {\bf 7}, it suffices to show that   $H^{1}(\X_{\text{\rm fppf}}, 
\Aut(A_{X}))$ is trivial. Since  $\Aut(A_{X})$ is smooth,  we may replace the 
fppf by the \'etale topology. Now apply
Proposition 5 (ii).  
\qed
\enddemo

\medskip

\noi {\bf Example 9.} If $A$ is the associative unital algebra of $n$ 
by $n$ matrices with entries in an algebraically closed field $k$, then all its covering 
algebras are  trivial. Indeed. $\Aut(A) \simeq \text{\rm \bf PGL}_{n}.$

\medskip
\proclaim{Proposition 10} Let $k$ be an algebraically closed  field of characteristic $0$. Assume 
that the $k$-algebra $A$ is finite 
dimensional,    and that the 
linear algebraic $k$-group $\AutA$ coincides with the group of 
automorphisms $\AutG$ of some semisimple algebraic $k$-group $G$ of 
either adjoint or simply connected type. Set 
$X = \Spec k[t,t^{-1}]$. Then there exists 
a canonical bijection between the following four sets:
\smallskip
\item{\text{\rm (1)}}  $H^{1}(X, \Aut(G_{X}))$

\item{\text{\rm (2)}} $H^{1}(X, \Out(G_{X}))$

\item{\text{\rm (3)}} Conjugacy classes of the (abstract) finite group 
$\text{\rm Out}(G) := \OutG(k).$

\item{\text{\rm (4)}} Isomorphism classes in $\Ralg$ of loop algebras 
of $A.$
\smallskip
\noi In particular, all these sets are finite. If in addition in the group $\text{\rm Out}(G)$  all  
elements are conjugate to their inverses, and if the centroid of all loop algebras of $A$ coincides 
with $R$, then (4) above is equivalent to

\item{\text{\rm (5)}} {\it Isomorphism classes in $k$-alg of loop algebras 
of $A.$}

\endproclaim
\demo{Proof}
\noi (1) $\simeq$ (2). This was established in Theorem 3 (iv).

\noi (2) $\simeq$ (3). Let $E$ be the set of continuous group homomorphisms from 
the fundamental group $\Pi(X,x)$ into $\text{\rm Out}(G)$. Clearly $\text{\rm Out}(G)$ acts on 
$E$ by conjugation. Since $\Out(G_{X})$ is finite and constant, $H^{1}(X, 
\Out(G_{X}))$ can be 
computed as the quotient set $E/\text{\rm Out}(G)$  [SGA1].  Now from 
Proposition 4 we 
know that $\Pi(X,x) \simeq \text{inv lim}\, \Bbb Z/n\Bbb Z$. It follows 
that the elements of $E$ can be identified in a natural way with elements of $\text{\rm 
Out}(G)$, and that then two elements of $E$ are equivalent under the 
action of $\text{\rm Out}(G)$ if and only if the corresponding 
elements of $\text{\rm Out}(G)$ are conjugate.

(3) $\simeq$ (4).  First note that since we are in characteristic $0$ 
all loop algebras are of the form $L(\s)$ for some $\s \in \text{\rm 
Aut}(A) = \Aut(A)(k) = \AutG(k)$. The correspondence assigns 
to the $R$-isomorphism class of $L(\s)$ the conjugacy class of 
${\overline\s}$. That such map exists and is bijective follows from the interpretation of the 
$L(\s)$'s as torsors, the explicit construction of $L^{1}(\s)$ 
given in {\bf 7}, and Proposition 3(iv).

\medskip
Assume now that $A$ above is such that in the group $\text{\rm Out}(G)$  all  
elements are conjugate to their inverses, and also such that the 
centroid $Z$ of all the 
$L(\s)$'s (viewed as $k$-algebras) coincide with $R$ (namely 
$Z = k[z^{m},z^{-m}]  \simeq R$ where $m$ is a chosen period of $\s$).   
We must show that  if two loop algebras $L(\s{_1})$ and $L(\s{_2})$ 
are isomorphic as $k$-algebras then they are isomorphic as 
$R$-algebras.  Since 
all the isomorphism classes in question do not depend on the choice of period, 
 we 
may assume that both $\s{_1}$ and $\s{_2}$ have period $m$.
We now reason as in [ABP]. Fix a $k$-algebra isomorphism $\psi$ between $L(\s{_1})$ and 
$L(\s{_2})$. Then $\psi$ induces 
an automorphism $\psi_{Z}$ of the centroid $Z$ of $L(\s{_1})$. 
There are two ``types'' of automorphisms of $Z \simeq R$: those that 
do not interchange $kt$ and $kt^{-1}$, and those which do . The upshot 
of this is that $L(\s{_1})$   is isomorphic as an $R$-algebra to either 
$L(\s_{2})$ 
or $L(\s_{2}^{-1})$ depending on the type of $\psi_{Z}$. 

Having establish this, we appeal to (3) $\simeq$ (4) to conclude that ${\overline\s_{1}}$ and 
${\overline\s_{2}}$ are conjugate (because in $\text{\rm 
Out}(G)$ elements are conjugate to their inverses), and hence that $L(\s{_1})$ and $L(\s{_2})$ 
are isomorphic as $R$-algebras (again by (3) $\simeq$ (4)).
\qed

\enddemo

We can now recover the  parts  of Theorem 2 concerning loop 
algebras.
\medskip
\proclaim{Corollary 11} Let $\g$ be a finite  dimensional simple Lie 
algebra over an algebraically closed field $k$ of characteristic $0.$ 
Let $\s_{1}$ 
and $\s_{2}$ be two 
automorphisms  of $\g$ of finite order. For $\C(\s_{1})$ to be isomorphic 
to $\C(\s_{2})$ as 
algebras over $k$, it is necessary and sufficient that 
$\overline  \s_{1}$ and $\overline \s_{2}$ be conjugate in $\text{\rm Out}(\g).$
\endproclaim

\demo{Proof} Let $G$ be the simply connected Chevalley-Demazure group  
corresponding to $\g.$ Then  $\Aut(G_{X}) \simeq \Aut(\g_{X})$ with 
$\Out(G)$ corresponding to $\Out(\g)$ ([SGA3]  XXV). Furthermore 
$\text{\rm Out}(G)$ is the group of automorphisms of the 
corresponding Dynkin diagram, and one knows that in these groups all 
elements are conjugate to their inverses. That  the centroids 
of the loop algebras of $\fg$ coincide with $R$  
is easy to check (see [ABP] for details).  \qed
\enddemo

\noi {\bf Remark 12} It is natural to ask if Proposition 10 holds for 
symmetrizable Kac-Moody Lie algebras. We look at this problem in 
[ABP], and  make good progress  by using  the Gantmacher-like decomposition of 
automorphisms  described in [KW].  

It would be interesting to know if the answer to this question can be 
had by 
purely cohomological methods (as the finite dimensional case here as 
well as [Pzl1] and 
[Pzl2] seem to suggest is possible). 
The abstract construction of section {\bf 7} applies, but the real difficulty 
of course comes when one tries to recreate Proposition 3 for the 
various group 
schemes attached to Kac-Moody algebras (See [Tts]). This appears to be  a very difficult 
question, but in the affine case at least, progress should be possible.

\medskip
\noi {\bf Acknowledgement:} I am thankful  to Philippe 
Gille for his many useful suggestions and  comments. 

I  also wish to thank the people at the Fields Institute for their 
hospitality while most of this work was done, and David Harari for 
his inspiring six lectures on torsors.

\end